\documentclass[12pt]{amsart}
\usepackage{amsthm, amsfonts, amssymb, amsmath, amscd}
\usepackage{graphicx}
\usepackage{mathptmx}
\usepackage[margin=3.2cm]{geometry}
\numberwithin{equation}{section}

\newcommand{\X}{\mathbf{X}}
\newcommand{\Y}{\mathbf{Y}}
\newcommand{\Z}{\mathbb{Z}}

\newcommand{\Q}{\mathbb{Q}}
\newcommand{\R}{\mathbb{R}}

\newcommand{\T}{\mathrm{(T)}}

\newcommand{\SL}{\mathrm{SL}}

\newcommand{\FX}{\mathrm{FX}}

\newcommand{\To}{\longrightarrow}
\newcommand{\e}{\varepsilon}
\newcommand{\Mod}{\mathrm{Mod}}
\newcommand{\cat}{\textsc{CAT(0)}}

\newtheorem{thm}{Theorem}[section]

\newtheorem{lem}[thm]{Lemma}
\newtheorem{cor}[thm]{Corollary}
\newtheorem{conj}[thm]{Conjecture}
\newtheorem{prop}[thm]{Proposition}

\theoremstyle{definition}
\newtheorem{defn}[thm]{Definition}
\newtheorem{rem}[thm]{Remark}

\newtheorem{ex}[thm]{Example}

\newtheorem{prob}[thm]{Problem}

\begin{document}
\title{Group actions on median spaces}
\author{Bogdan Nica}
\address{Department of Mathematics, Vanderbilt University, Nashville, TN 37240}
\email{bogdan.nica@vanderbilt.edu}
\keywords{}
\date{\today}

\begin{abstract}
We investigate the geometry of median metric
spaces. The group-theoretic applications are towards Kazhdan's property (T) and Haagerup's property.
\end{abstract}

\maketitle

\section{Introduction}
We are interested in Kazhdan's property $\T$ and its opposite, Haagerup's property, for discrete countable groups. These two properties can be defined in many ways, but the relevant perspective for this paper is that of isometric actions on (real) Hilbert spaces. A group $G$ is a \emph{Kazhdan group} if every isometric action of $G$ on a Hilbert space has a fixed point. A group $G$ is a \emph{Haagerup group} if $G$ admits a proper isometric action on a Hilbert space. These two properties are opposite in the sense that the groups which are both Kazhdan and Haagerup are precisely the finite groups. 

The Kazhdan vs. Haagerup dichotomy arises by fixing a geometry, the Hilbert space geometry in this case. On one side we isolate groups which are incompatible with the geometry, i.e., having only trivial actions. On the other side we have the groups which manifest their compatibility with the given  geometry by admitting proper actions. The archetype can be formulated as follows. Fix a family $\X$ of (isometry classes of) metric spaces, thought of as a geometry. Using $\X$ as a test-ground for isometric actions, we distinguish ``$\X$-rigid'' groups, respectively ``$\X$-amenable'' groups.

\begin{defn}\label{defFX}
A group has \emph{property FX} if each isometric action on a member of $\X$ is bounded.
\end{defn}

\begin{defn}
A group has \emph{property PX} if it admits a proper isometric action on a member of $\X$.
\end{defn}

\noindent Recall, an isometric action of a group $G$ on a metric space $X$ is said to be

$\cdot$ \emph{bounded} if, for some (all) $x\in X$, the orbit $\{gx:\: g\in G\}$  is bounded;

$\cdot$ \emph{proper} if, for some (all) $x\in X$, the ``coarse stabilizer'' $\{g\in G:\:|x\;gx|\leq R\}$ is finite for each $R\geq 0$.

The vestigial ``F'' in Definition~\ref{defFX} suggests a globally fixed point. In certain geometries, one can pass from bounded orbits to fixed points by showing that bounded sets have unique ``centers''. This is indeed the case for complete, uniformly convex metric spaces (see Appendix~\ref{UC}).

Property FX and property PX are opposite in the same sense as before, namely the groups enjoying both properties are precisely the finite groups. Note that property FX is inherited by quotients, while property PX is inherited by subgroups. In particular, an extension $1\rightarrow F \rightarrow G \rightarrow P\rightarrow 1$, where $F$ is an infinite FX group and $P$ is an infinite PX group, yields an ``$\X$-ambiguous'' group $G$ lacking both FX and PX. 

The better known avatars of property FX are property FA, relative to the family of simplicial trees, and property FH, relative to the family of Hilbert spaces. As for property PX, it has received most attention for the family Hilbert spaces. The variations are numerous. For instance, one can consider the F$L_p$ / P$L_p$ properties for $1\leq p< \infty$. The chosen geometry can be that of CAT(0) metric spaces, or maybe that of CAT(0) cube complexes; that of uniformly convex Banach spaces, or that of ultrametric spaces, et cetera.

The basic connection between action on trees  (FA vs. PA) and actions on Hilbert spaces (Kazhdan vs. Haagerup) inspires much of the ideas in the FX vs. PX framework. To gain understanding, one often considers a second geometry  $\Y$ and tries to relate property FX, respectively PX, to property FY, respectively PY. Such relationships, typically established at the metric level between members of $\X$ and members of $\Y$, allow for group-theoretic insights.

In this paper, we consider the FX vs. PX phenomenon relative to the median geometry. More precisely, we consider the family of median spaces as well as the subfamily of median graphs. We relate to the Kazhdan and the Haagerup properties as follows (Propositions~\ref{FH_FM} and ~\ref{PM_PH}):

\begin{thm} Every isometric action of a Kazhdan group on a median space is bounded. A group that admits a proper action on a median space has the Haagerup property.
\end{thm}

This paper is a revision of work submitted in 2004 as a Master's thesis at McGill University (\cite{Nic}). We have expanded and updated the introductory sections, which present a very selective and incomplete survey of FX-like and PX-like properties. The part on median spaces (Section 5 here) has been trimmed and slightly rewritten. Previous work of Verheul is acknowledged in this section. The part devoted to the correspondence between spaces with walls and CAT(0) cube complexes has already been published (\cite{Nic04}). To avoid repetition, we have excised the material on ``cubulations'', so the discussion around median graphs now makes a rather concise and expository Section 6. Finally, we have added an Appendix. The first item in the Appendix is a general Circumcenter Lemma for uniformly convex metric spaces; this might be known, but we have not been able to trace it in the literature. The second item is a proof from The Book for the Mazur-Ulam theorem.

Conjecture~\ref{conject} was recently settled by Chatterji, Dru\c{t}u and Haglund in \cite{CDH07}. In light of the renewed interest in median geometry, we have decided to make this work available.

\textbf{Acknowledgments:} I am grateful to Dani Wise for being an inspiring adviser, and to Mark Sapir and Guoliang Yu for encouraging me to publish this paper. During 2002-2004, I was supported by NSERC (National Sciences and Engineering Council of Canada).

\section{Property FX}
This section surveys some aspects of the FX landscape. We consider groups that are rigid with respect to one the following geometries:

$\cdot$ trees and $\R$-trees (properties FA and FA$_\R$)

$\cdot$ CAT(0) complexes of fixed dimension (property FA$_n$)

$\cdot$ Hilbert space geometry (Kazhdan's property (T), or property FH)

$\cdot$ $L_p$ geometry for $1\leq p<\infty$ (property FL$_p$)

\noindent The review of Kazhdan groups is inevitably brief. For more on property (T), see \cite{HV89} and its follow-up, \cite{BHV}. Shalom's survey \cite{Sha06} is also recommended.

\subsection{Property FA and Property FA$_\R$} Serre~\cite{Ser80} introduced what is perhaps the earliest avatar of property $\FX$, namely geometric rigidity with respect to simplicial trees:

\begin{defn} A group $G$ has \emph{property FA} if every isometric action of $G$ on a
simplicial tree fixes a vertex or an edge.
\end{defn}

Considering $\R$-trees (``arbres r\'{e}els'') instead of simplicial trees (``arbres''), one has the following analogue of property FA:

\begin{defn}\label{far} A group $G$ has \emph{property FA$_\R$} if every isometric action of $G$ on an $\R$-tree has a fixed point.
\end{defn}

The latter property has been considered by Culler and Vogtmann \cite{CV96}, and by Bridson \cite{Bri08}. Both \cite{CV96} and \cite{Bri08} call it ``property F$\R$'' but, under our FX notation, that would mean geometric rigidity with respect to the real line. For consistency's sake, we henceforth use the terminology of Definition~\ref{far}.

Every simplicial tree can be turned into an $\R$-tree, by making each edge isometric to the unit interval. Furthermore, isometric actions extend uniquely. The following is then clear:

\begin{prop} Property FA$_\R$ implies property FA. 
\end{prop}

It is not known whether the converse holds. In an equivalent form, this question seems to have been first raised by Shalen \cite[Question A, p.286]{Sha87}. 

As shown by Serre, property FA has a purely algebraic reformulation:

\begin{thm}\label{algebraic_FA}
A group $G$ has property FA if and only $G$ is finitely generated, $G$ has finite abelianization, and $G$ is not a non-trivial amalgamated product.
\end{thm}

A group is said to be \emph{indicable} if it admits an epimorphism to $\Z$. The most immediate examples of finitely generated groups that lack property FA are to be found among indicable groups. For example, a finite presentation with more generators than relators defines an  indicable group. 

\begin{prop}\label{socioFA} Property FA$_\R$ and property FA are preserved by taking quotients, extensions, and finite index supergroups. 
\end{prop}

That property FA$_\R$ and property FA are not inherited by finite index subgroups will become apparent from the examples below.

\begin{ex}\label{ExFAR} Here are some groups that enjoy property FA$_\R$:
\begin{enumerate}
\item Finitely generated torsion groups

\item Coxeter groups $\langle s_1, \dots, s_n | (s_is_j)^{m_{ij}}=1\rangle$ with all $m_{ij}$'s finite

\item $\textrm{SL}_n(\Z)$ for $n\geq 3$, together with its finite index subgroups 

\item $\textrm{E}_n(A)$, the elementary linear group of a finitely generated ring $A$, for $n\geq 3$

\item $\textrm{Aut}(F_n)$ for $n\geq 3$

\item $\Mod(S_g)$, the mapping class group of a closed orientable surface, for $g\geq 2$
\end{enumerate}

\noindent Items i) \& ii) are due to Serre \cite[Cor. 2, p.64]{Ser80}; they are stated for property FA, but the approach applies to $\R$-trees as well. Items iii) to vi) are due to Culler and Vogtmann \cite{CV96}. For property FA, iii) is due to Serre \cite[\S I.6.6]{Ser80}, and v) is due to Bogopolski \cite{Bog87}.
See also Bridson \cite{Bri08} for a different approach to iii) and v).

All but one of the exceptions in ii)-vi) are easy to explain. In ii), if some $m_{ij}=\infty$ then the Coxeter group splits as a non-trivial amalgam. In iii), $\textrm{SL}_2(\Z)$ has a well-known amalgam decomposition; in vi), $\Mod(S_1)$, the mapping class group of the torus, is isomorphic to $\textrm{SL}_2(\Z)$. In v), $\textrm{Aut}(F_2)$ maps onto $\textrm{GL}_2(\Z)$, which does not have property FA. As for the exception in iv), namely $\textrm{E}_2(A)$, it is rather elusive. The manageable case is that when $\textrm{E}_2(A)=\textrm{SL}_2(A)$, for much is known about the group $\textrm{SL}_2(A)$ for various familiar rings $A$. In this case, property FA for $\textrm{E}_2(A)$ may or may not hold; for instance, it holds when $A$ is the ring of integers in $\Q(\sqrt{-3})$ or in $\Q(\sqrt{3})$, and it fails when $A$ is $\mathbb{F}_q[X]$, or $\Z[1/p]$, or the Gaussian integers $\Z[i]$.

As property FA$_\R$ is inherited by quotients, and $\mathrm{E}_n(A)\twoheadrightarrow \mathrm{E}_n(B)$ whenever $A\twoheadrightarrow B$, we have that iv) is equivalent to the apparently weaker statement that $\textrm{E}_n(\Z[X_1,\dots, X_d])$ has property FA$_\R$ for $n\geq 3$ and $d\geq 0$. Via Suslin's $\SL_n(\Z[X_1,\dots, X_d])=\textrm{E}_n(\Z[X_1,\dots, X_d])$, this is the same as saying that the ``universal lattice'' $\SL_n(\Z[X_1,\dots, X_d])$ has property FA$_\R$ for $n\geq 3$ and $d\geq 0$.

Some of the groups listed above have finite index subgroups that map onto $\Z$. Every infinite finitely generated Coxeter group is virtually indicable (\cite{Gon97}, \cite{CLR98}). McCool \cite{Mc89} showed that $\textrm{Aut}(F_3)$ is virtually indicable. It is an open question whether finite index subgroups of $\textrm{Aut}(F_{n\geq 4})$ have property FA$_\R$.
\end{ex}

\subsection{Property FA$_n$} There are several ways of defining higher analogues of property FA. The following version is due to Farb \cite{Fa08}.
\begin{defn} A group $G$ has \emph{property FA$_n$} if every isometric action of $G$ on a
complete, CAT(0) $n$-complex has a fixed point.
\end{defn}

An \emph{$n$-complex}, in this definition, is an $n$-dimensional cell complex with finitely many isometry types of cells; an $n$-complex is not required to be locally finite. 

Note that one recovers the usual property FA in the case $n=1$. Note also that property FA$_n$ implies property FA$_{n-1}$. As in Proposition~\ref{socioFA}, property FA$_n$ is preserved by taking quotients, extensions and finite-index supergroups.

\begin{ex}\label{exFAn}
Interesting examples from \cite{Fa08} include:
\begin{enumerate}
\item $\SL_n(\Z[X_1,\dots,X_d])$ has property FA$_{n-2}$ for $n\geq 3$ and $d\geq 0$
\item $\SL_{n+1}(\Z[1/p])$, and groups generated by reflections in the sides of a euclidean or hyperbolic $n$-simplex of finite volume, satisfy FA$_{n-1}$ but fail FA$_n$
\end{enumerate}
\end{ex}

\subsection{Property FH} Rigidity of actions relative to (real) Hilbert space geometry leads to Serre's property FH. By the Delorme-Guichardet theorem, property FH is equivalent to Kazhdan's property (T).  We may then use the following terminology:

\begin{defn}
A group $G$ has is a \emph{Kazhdan group} if every isometric action of $G$ on a Hilbert space has a fixed point.
\end{defn}

Note that an isometric action of a group on a real normed space is automatically affine, by the Mazur-Ulam theorem (see Appendix~\ref{MU}).

\begin{prop} 
The collection of Kazhdan groups is closed under taking quotients, extensions, and finite-index subgroups.
\end{prop}

The connection between property FH and property FA is due to Watatani \cite{Wat82}:

\begin{thm}\label{kazhdan_FA}
Kazhdan groups have property FA.
\end{thm}

In particular, Kazhdan groups are finitely generated, and have finite abelianization. This was first proved in \cite{Kaz67} in a way that does not go through property FA. In fact, these two consequences of property (T) motivated Kazhdan in introducing his fundamental property.

Noskov \cite{Nos93} proved the following result, presumably a generalization of Theorem~\ref{kazhdan_FA}:

\begin{thm}\label{kazhdan_FAR}
Kazhdan groups have property FA$_\R$.
\end{thm}

\begin{ex} The following groups are listed in Example~\ref{ExFAR} as having property FA$_\R$.
\begin{enumerate}
\item Finitely generated torsion groups

\item Coxeter groups $\langle s_1, \dots, s_n | (s_is_j)^{m_{ij}}=1\rangle$ with all $m_{ij}$'s finite

\item $\textrm{SL}_n(\Z)$ for $n\geq 3$, together with its finite index subgroups 

\item $\textrm{E}_n(A)$, the elementary linear group of a finitely generated ring $A$, for $n\geq 3$

\item $\textrm{Aut}(F_n)$ for $n\geq 3$

\item $\Mod(S_g)$, the mapping class group of a closed orientable surface, for $g\geq 2$
\end{enumerate}
Which ones are Kazhdan groups?

i) There are three significant sources of infinite, finitely generated torsion groups. Firstly, Grigorchuk-type groups acting on rooted trees. In this family, most groups are amenable, but it is still not known whether a Kazhdan group can exist. Secondly, Golod-Shafarevich groups. There are Kazhdan groups in this family, as shown by Ershov \cite{Ersh}. On the other hand, it is not known whether a Golod-Shafarevich group can be amenable. Thirdly, free Burnside groups $B_n^m=\langle x_1,\dots, x_n|\; w^m=1\textrm{ for all }w=w(x_1,\dots,x_n)\rangle$. It is an open question whether some $B_n^m$ is a Kazhdan group. One should mention here that $B_n^m$ is non-amenable for $n\geq 2$ and $m$ sufficiently large (Novikov-Adian, Olshanski, Ivanov, Lysionok).

ii) No infinite Coxeter group is a Kazhdan group; this is due to Bo\.{z}ejko, Januszkiewicz and Spatzier \cite{BJS88}. In fact, their proof shows that Coxeter groups are Haagerup groups.

iii)  For $n\geq 3$, $\textrm{SL}_n(\Z)$ is a Kazhdan group. This fundamental example is due to Kazhdan \cite{Kaz67}. Kazhdan's proof views $\textrm{SL}_n(\Z)$ as a lattice in $\textrm{SL}_n(\R)$. Shalom's proof in \cite{Sha99} does not use the ambient Lie group; instead, it uses the bounded generation of $\textrm{SL}_{n\geq 3}(\Z)$ by elementary matrices.

iv) $\textrm{E}_n(A)$ is a Kazhdan group for $n\geq 3$. Shalom \cite{Sha06} proved that $\textrm{E}_n(A)$ is a Kazhdan group for $n\geq 2+\dim(A)$, where $\dim(A)$ denotes the Krull dimension of $A$. Vaserstein (unpublished) improved Shalom's result to the present form. 

v) It is an important open problem whether $\mathrm{Aut}(F_n)$ is a Kazhdan group for $n\geq 4$. Recall that $\mathrm{Aut}(F_3)$ is virtually indicable, so it cannot be a Kazhdan group.

vi) Andersen \cite{And07} has recently shown that $\Mod(S_g)$ is not a Kazhdan group for $g\geq 2$.
\end{ex}

\begin{rem} 
Note that a Kazhdan group need not have Farb's property FA$_n$ for $n\geq 2$. This is witnessed by $\SL_{n+1}(\Z[1/p])$: on one hand it is a Kazhdan group, on the other hand it fails property FA$_n$ (Example~\ref{exFAn}\:ii).
\end{rem}

\subsection{Property FL$_p$ ($1\leq p<\infty$)} Along with another $L_p$-version of Kazhdan's property (T), this property is considered by Bader, Furman, Gelander and Monod \cite{BFGM}. They prove, among many other things, the following:

\begin{thm} Let $G$ be a Kazhdan group. Then $G$ has property FL$_p$ for $1\leq p< 2+\e$, where $\e=\e(G)>0$.
\end{thm}

The $p=1$ case was previously proved by Robertson and Steger. The $2\leq p< 2+\e$ part is unpublished work of Fisher and Margulis.

The next result is a particular case of \cite[Thm.B]{BFGM}:

\begin{thm} For $n\geq 3$, $\SL_n(\Z)$ has property FL$_p$ for all $1\leq p<\infty$.
\end{thm}

\section{Property PX}
We now look at several manifestations of property PX. The geometries we consider are the following:

$\cdot$ Hilbert space geometry (Haagerup's property, or a-T-menability)

$\cdot$ CAT(0) cube complexes (property PQ)

$\cdot$ $L_p$ geometry for $1\leq p<\infty$ (property PL$_p$)

\noindent For more on Haagerup groups, see \cite{CCJJV}.

\subsection{Property PH} The property of having a proper isometric action on a Hilbert space was explicitly formulated by Gromov \cite{Gro93}, who called it \emph{a-T-menability}. It turned out that a-T-menability is equivalent to an analytic property for groups that had been previously considered, the \emph{Haagerup approximation property}. Both this property and the Rapid Decay property (which does not concern us here) stem from Haagerup's paper \cite{Haa79}. 

We adopt the following:

\begin{defn}
A group $G$ is a \emph{Haagerup group} if $G$ admits a proper isometric action on a Hilbert space.
\end{defn}

\begin{ex}\label{amenable}
Amenable groups are Haagerup groups. 

Here is a proof taken from \cite{BCV95}. Let $G=\{g_1,g_2, \ldots\}$ be an amenable group. By F$\o$lner's characterization, there is a sequence $(F_n)_{n\geq 1}$ of finite subsets of $G$ such that
\begin{displaymath}
\frac{|gF_n\triangle F_n|}{|F_n|}\To 0
\end{displaymath}
for every $g\in G$. We need, however, a better asymptotic control over the F$\o$lner sequence
$(F_n)_{n\geq 1}$. By passing to a subsequence if necessary, we
may assume that for each $n\geq 1$ we have the following:
\begin{displaymath}
\frac{|g_iF_n\triangle F_n|}{|F_n|}<\frac{1}{n^3}\quad \textrm{ for }  1\leq i\leq n
\end{displaymath}
For each $n\geq 1$, turn the usual linear isometric action of $G$ on
$\ell_2(G)$ into an affine one
\begin{displaymath}
g*_n \phi=g \phi+\sqrt{\frac{n}{|F_n|}}(\chi_{_{gF_n}}-\chi_{_{F_n}})
\end{displaymath}
and bundle all these isometric actions into a single isometric action on the $\ell_2$-sum 
$\bigoplus_n\ell_2(G)$ in the obvious way:
\begin{displaymath}
g* (\phi_n)_{n\geq 1}=(g*_n \phi_n)_{n\geq 1}
\end{displaymath}
This action is well-defined, since for each $g=g_i$ we have:
\begin{eqnarray*}
\sum_{n\geq 1}\bigg\|\sqrt{\frac{n}{|F_n|}}(\chi_{_{g_iF_n}}-\chi_{_{F_n}})\bigg\|^2
&=& \sum_{n\geq 1}\frac{n}{|F_n|}|g_iF_n\triangle F_n|
= \sum_{n<i}n\;\frac{|g_iF_n\triangle F_n|}{|F_n|}+\sum_{n\geq
i}n\;\frac{|g_iF_n\triangle F_n|}{|F_n|}\\
&<&\sum_{n<i}n\;\frac{|g_iF_n\triangle F_n|}{|F_n|}+\sum_{n\geq i}\frac{1}{n^2}<\infty
\end{eqnarray*}
Furthermore, the action is proper. Indeed, let $N$ be a fixed positive integer. There are only finitely many $g\in G$ such that $gF_N$ meets $F_N$. For any other $g$ we have $|gF_N\triangle F_N|=2|F_N|$, hence the following estimate:
\begin{eqnarray*}
\|g*0\|^2=\sum_{n\geq 1}\bigg\|\sqrt{\frac{n}{|F_n|}}(\chi_{_{gF_n}}-\chi_{_{F_n}})\bigg\|^2=\sum_{n\geq
1}\frac{n}{|F_n|}|gF_n\triangle F_n| \geq \frac{N}{|F_N|}|gF_N\triangle F_N|=2N
\end{eqnarray*}
We conclude that, for each $N$, the set $\{g\in G
:\|g*0\|<\sqrt{2N}\}$ is finite. \end{ex}

\begin{prop}
The collection of Haagerup groups is closed under taking subgroups, extensions with amenable quotients, and ascending countable unions.
\end{prop}

The family of Haagerup groups is not closed under extensions. For example, $\Z^2$ and $\SL_2(\Z)$ are Haagerup groups yet $\Z^2 \rtimes \SL_2(\Z)$ is not a Haagerup group. The reason is that every isometric action of $\Z^2 \rtimes \SL_2(\Z)$ on a Hilbert space has a point fixed by $\Z^2$, i.e., $\Z^2\rtimes \SL_2(\Z)$ is Kazhdan relative to $\Z^2$. Note, however, that $\Z^2 \rtimes \SL_2(\Z)$ is not Kazhdan. 

A rich supply of Haagerup groups is given by groups acting properly on CAT(0) cube complexes. This is discussed in the next item.

\subsection{Property PQ} In the following definition, the letter Q is meant to suggest ``cube'':

\begin{defn} A group $G$ has \emph{property PQ} if $G$ has a proper action on a CAT(0) cube complex.
\end{defn}

The roots of the next result lie in \cite{NRe97}:

\begin{thm}\label{PQ-PH} If $G$ has property PQ, then $G$ is a Haagerup group.
\end{thm}

\begin{ex} The groups enjoying property PQ include:
\begin{enumerate}
\item virtually free, finitely generated groups
\item lattices in products of regular trees
\item Coxeter groups (\cite{NRe03})
\item the Thompson group (\cite{Far03})
\item $C'(1/6)$ small cancellation groups (\cite{Wis04})
\end{enumerate}
\end{ex}
It should be noted that free groups are the original examples of Haagerup groups (\cite{Haa79}). That Coxeter groups are Haagerup is already implicit in \cite{BJS88}.

Examples of Haagerup groups that do not have property PQ were given recently by Haglund \cite{Hag07}. He shows that finitely generated groups containing a distorted infinite cyclic subgroup cannot act properly on a CAT(0) cube complex. On the other hand, Haagerup groups with distorted copies of $\Z$ are easy to find. The discrete Heisenberg group
\[\mathrm{UT}_3(\Z)=\Bigg\{\begin{pmatrix}
1 & x & y\\
& 1 & z\\
& & 1
\end{pmatrix}: x,y,z\in \Z\Bigg\}
\]
is amenable, hence Haagerup, and $\langle y\rangle $ is a distorted subgroup. The Baumslag-Solitar group 
\[BS_{p,q}=\langle a,b| \: ab^p a^{-1}=b^q\rangle\]
is a Haagerup group (Gal-Januszkiewicz), and the subgroup $\langle b \rangle$ is distorted when $p\neq q$.

\subsection{Property PL$_p$ ($1\leq p<\infty$)}
An immediate adaptation of Example~\ref{amenable} shows that amenable groups have property PL$_p$ for each $p\geq 1$. Also, it is not hard to see that groups acting properly on CAT(0) cube complexes have property PL$_p$ for each $p\geq 1$. However, more is true, as observed by de Cornulier, Tessera and Valette \cite[Prop. 3.1]{CTV06}. The proof of the following fact relies on a previously known characterization of Haagerup property in terms of measured spaces with walls.

\begin{thm} Let $G$ be a Haagerup group. Then $G$ has property PL$_p$ for all $p\geq 1$.
\end{thm}

More interesting is the following result, due to Yu \cite{Yu05}:

\begin{thm}\label{yu} Let $G$ be a hyperbolic group. Then $G$ has property PL$_p$ for $p$ sufficiently large.
\end{thm}

It should be noted that many hyperbolic groups are Kazhdan groups, e.g., cocompact lattices in $\mathrm{Sp}(n,1)$ for $n\geq 2$. According to \cite{CTV06}, for such lattices one can take $p>4n+2$ in Theorem~\ref{yu}.

\begin{rem} For some geometries $\X$, the corresponding properties PX and FX are trivial.

If $\X$ consists of bounded metric spaces, then every group has property FX, and the groups enjoying property PX are precisely the finite groups. On the other hand, if $\X$ is a ``large'' geometry then the situation is reversed: every group has property PX and, consequently, the groups enjoying property FX are precisely the finite groups. This is the case when $\X$ consists of the (separable, non-reflexive) space $c_0$, or the (non-separable, non-reflexive) space $\ell_\infty$. It is also the case when $\X$ consists of the strictly convex Banach spaces.

Indeed, let $L$ be a proper length function on $G$ (any countable group inherits a length function from an embedding into a finitely generated group). Modify the regular isometric action $(g,\phi)\mapsto g\phi$ of $G$ on $c_0(G)$ into the following isometric action of $G$ on $c_0(G)$:
\[g*\phi=g\phi +g\sqrt{L}-\sqrt{L}\]
To see that this action is well-defined, for each $g\in G$ write 
\[g\sqrt{L}-\sqrt{L}=\frac{gL-L}{g\sqrt{L}+\sqrt{L}}\]
and note that the numerator is bounded (by $L(g)$) whereas the denominator is a proper function on $G$. Thus $g\sqrt{L}-\sqrt{L}\in c_0(G)$. Furthermore, $\|g*0\|_\infty=\sqrt{L(g)}$ so the above action is proper. This shows that every group has a proper isometric action on $c_0$.

The same action as above, but this time interpreted on $\ell_\infty (G)$, shows that every group has a proper isometric action on $\ell_\infty$.

Brown and Guentner \cite{BG05} show that every group admits a proper isometric action on some strictly convex Banach space. The proof is an elaboration of the strategy used in Example~\ref{amenable}.
\end{rem}


\section{Median algebras}
The purpose of this short section is to define median algebras and to present some of their key features. There are (at least) two equivalent ways of defining a median algebra. The algebraic perspective is to define median algebras as sets endowed with a ternary operation, called the median operation, that satisfies certain axioms. Better adapted to our geometric context is the definition we adopt below:

\begin{defn}\label{MA}
A \emph{median algebra} is a set $X$ equipped with an interval map 
$[\cdot,\cdot]:X\times X\To\mathcal{P}(X)$ such that for all $x,y,z\in X$ the following are satisfied:

($\textrm{MA}_1$)\quad $[x,x]=\{x\}$

($\textrm{MA}_2$)\quad $[x,y]=[y,x]$

($\textrm{MA}_3$)\quad if $z\in [x,y]$ then $[x,z]\subseteq [x,y]$

($\textrm{MA}_4$)\quad $[x,y]$, $[y,z]$, $[z,x]$ have a unique common point, called the \emph{median} of $x,y,z$

\noindent 
A map $f:X\To X'$ between median algebras is a \emph{median morphism} if it respects the interval structure, in the sense that $f\big([x,y]\big)\subseteq \big[f(x),f(y)\big]$ for all $x,y\in X$.
\end{defn}

An interval structure induces a notion of convexity:

\begin{defn} 
Let $X$ be a median algebra. A subset $A\subseteq X$ is \emph{convex} if $[x,y]\subseteq A$ for
all $x,y\in A$. A subset $H\subseteq X$ is a \emph{halfspace} if both $H$ and $H^c$ are convex.
\end{defn}

A crucial feature of halfspaces in a median algebra is the following separation property. For a proof, see \cite[\S 2]{Rol98}. 

\begin{thm}\label{separate}
Let $X$ be a median algebra and $C_1, C_2\subseteq X$ be disjoint convex sets. Then there is a halfspace $H\subseteq X$ separating $C_1$ and $C_2$, i.e.,  $C_1\subseteq H$ and $C_2\subseteq H^c$.
\end{thm}

\begin{ex}[Boolean median algebra]\label{boolean_median}
Any power set $\mathcal{P}(X)$ is a median algebra under the interval assignment $[A,B]:=\{C: A\cap B\subseteq C\subseteq A\cup B\}$. The boolean median of $A$, $B$, $C$ is $(A\cap B)\cup(B\cap C)\cup(C\cap A)=(A\cup B)\cap(B\cup C)\cap(C\cup A)$. The non-empty halfspaces containing $X$ are precisely the ultrafilters on $X$. Recall, an \emph{ultrafilter} on $X$ is a collection $\mu$ of subsets of $X$ which satisfies the following properties: ($\textrm{U}_1$) $X\in\mu$, ($\textrm{U}_2$) $A,B\in\mu$ implies $A\cap B\in \mu$, ($\textrm{U}_3$) for all $A\subseteq X$, either $A\in\mu$ or $A^c\in\mu$.
\end{ex}

Any median algebra is isomorphic to a subalgebra of a boolean median algebra. Indeed, let $X$ be a median algebra, let $\mathcal{H}$ be the collection of halfspaces of $X$, and denote by $\sigma_x$
the collection of halfspaces containing $x\in X$. We obtain a map $\sigma:X\To\mathcal{P}(\mathcal{H})$ that is easily checked to be an injective median morphism.

We thus have a ``boolean method'' for proving (non-existential) statements about median algebras. The following result is an illustration of this method.

\begin{lem}\label{finite_median}
In a median algebra, the median closure of a finite set is finite.
\end{lem}

\noindent Given a median algebra $X$ and a subset $A\subseteq X$, the median closure of $A$ is the smallest subset of $X$ containing $A$ that is stable under taking medians.

\begin{proof} It suffices to prove the statement for boolean median algebras. A finite $\mathcal{A}\subseteq\mathcal{P}(X)$ generates a finite sublattice $\langle\mathcal{A}\rangle$ in $\big(\mathcal{P}(X), \cup,\cap\big)$; indeed, $\langle\mathcal{A}\rangle$ is obtained by taking unions of intersections of sets from $\mathcal{A}$. As $\langle\mathcal{A}\rangle$ is median stable, we are done.
\end{proof}

The following proposition gives an alternate definition of median morphisms. The relevance of this characterization, via halfspaces, becomes apparent in the context of spaces with walls (Definition~\ref{wallspacedef}). See \cite[Prop.3.7]{Nic04} for a proof.

\begin{prop}
Let $f:X\To X'$ be a map between median algebras. Then $f$ is a median morphism if and only if $f^{-1}(H')$ is a halfspace in $X$ whenever $H'$ is a halfspace in $X'$.
\end{prop}

\begin{rem}
a) The system of axioms given in Definition~\ref{MA} is a slight modification of a system of axioms due to Sholander. In Sholander's formulation, axioms ($\textrm{MA}_2$) and ($\textrm{MA}_3$) are replaced by the somewhat less natural

($\textrm{MA}_{2.5}$)\quad if $z\in [x,y]$ then $[x,z]\subseteq [y,x]$.

\noindent The definition of a median algebra in \cite{Nic04} has to be amended to Definition~\ref{MA} above. In \cite[Def.3.3]{Nic04}, the symmetry axiom ($\textrm{MA}_2$) is missing from the system ($\textrm{MA}_1$)--($\textrm{MA}_4$). An interval structure satisfying axioms ($\textrm{MA}_1$), ($\textrm{MA}_3$), ($\textrm{MA}_4$) but failing ($\textrm{MA}_2$) is, however, easy to find. Let $X=\{x,y,z\}$, and define the intervals of $X$ as follows:
\begin{align*}
[x,x]& =\{x\},\; [y,y]=\{y\},\; [z,z]=\{z\} & [y,z] & =[z,y]=\{y,z\}\\
[x,y]& =\{x,y,z\},\; [y,x]=\{x,y\} & [x,z]& =\{x,y,z\},\; [z,x]=\{x,z\}
\end{align*}

b) For more on median algebras, see Roller's study \cite{Rol98} and references therein.
\end{rem}

\section{Median spaces}
\subsection{Definition \& Examples} A metric space $X$ has a natural interval structure given by \emph{geodesic intervals}. Namely, $[x,y]:=\{t\in X: |x\;t|+|t\:y|=|x\;y|\}$ for $x,y\in X$. The geodesic interval structure of a metric space is easily seen to satisfy axioms ($\textrm{MA}_1$), ($\textrm{MA}_2$) and ($\textrm{MA}_3$) in the definition of a median algebra. Those metric spaces whose geodesic intervals satisfy axiom ($\textrm{MA}_4$) as well, thus giving them a median algebra structure, are the spaces we are interested in.

\begin{defn}
A metric space $X$ is \emph{median} if, for each triple $x,
y,z\in X$, the geodesic intervals $[x,y]$, $[y,z]$, $[z,x]$
have a unique common point.
\end{defn}

The following formula is often used: if $m\in X$ denotes the median of $x,y,z\in X$ then
\begin{equation}
|x\:m|=\frac{1}{2}\big(|x\:y|+|x\:z|-|y\:z|\big).
\label{leg}
\end{equation}
Note that the quantity on the right hand side is sometimes referred to as the Gromov inner product of $y$ and $z$ relative to $x$, denoted $<y,z>_x$.

If $X$ and $Y$ are median spaces, then
$X\times Y$ is a median space when equipped with the metric
$\big|(x_1,y_1)\;(x_2,y_2)\big|:=|x_1\:x_2|+|y_1\:y_2|$. This is easy to check.

The completion of a median space is median. We prove this in the next subsection.

\begin{ex}
$\R$-trees are median spaces. 
\end{ex}

\begin{ex} For any measure $\mu$, the space $L_1(\mu)$ of real-valued integrable functions is
median. 

Given $f,g\in L_1(\mu)$, then $m\in L_1(\mu)$ satisfies $\|f-m\|_1+\|m-g\|_1=\|f-g\|_1$ if and only if $|f(s)-m(s)|+|m(s)-g(s)|=|f(s)-g(s)|$ for $\mu$-almost all $s$. It follows that, for a given triple $f,g,h\in L_1(\mu)$, a function $m\in L_1(\mu)$ satisfies the geodesic conditions
\[\|f-m\|_1+\|m-g\|_1=\|f-g\|_1\]
\[\|g-m\|_1+\|m-h\|_1=\|g-h\|_1\]
\[\|h-m\|_1+\|m-f\|_1=\|h-f\|_1\]
if and only if $m(s)$ is the median of $f(s), g(s), h(s)$ for almost all $s$. That is, the median of $f,g,h\in L_1(\mu)$ is $m=(f\vee g)\wedge(g\vee h)\wedge(h\vee f)=(f\wedge g)\vee(g\wedge h)\vee(h\wedge f)$ $\mu$-a.e..
\end{ex}

\subsection{Completing median spaces}\label{CM} This section contains a proof of the fact that the completion of a median space is again a median space. 

We start with a general fact about median algebras:

\begin{lem}\label{colinear} Let $X$ be a median algebra, and let $m$ denote the median of $x,y,z\in X$. If $v\in [y,z]$ then $m\in[x,v]$.
\end{lem}

\begin{proof}\label{ineq} Assume that $m\notin[x,v]$. By Theorem~\ref{separate}, there is a halfspace $H$ such that $[x,v]\subseteq H$ and $m\in H^c$. If $y\in H$ then $m\in [x,y]\subseteq H$, a contradiction. Thus $y\in H^c$ and, similarly, $z\in H^c$. It follows that $v\in [y,z]\subseteq H^c$, contradicting the fact that $v\in H$.
\end{proof}

The next lemma holds the key to our desired result, Proposition~\ref{complete} below. Roughly speaking, the first inequality implies that a point which is close to fulfilling the role of median point, is in fact close to the median point. The second inequality yields the uniform continuity of the median function.

\begin{lem} Let $X$ be a median space.

i) If $m\in X$ is the median of $x,y,z\in X$, then for any $v\in X$ we have
\[|v\: m|\leq \big(|v\: x|+|v\: y|+|v\: z|\big)-\big(|m\: x|+|m\: y|+|m\: z|\big).\] 

ii) If $m\in X$ is the median of $x,y,z\in X$, and $m'\in X$ is the median of $x',y',z'\in X$, then 
\[|m\: m'|\leq |x\: x'|+|y\: y'|+|z\: z'|.\]
\end{lem}

\begin{proof}
i) Let $n$ be the median of $v,x,y$. Then:
\begin{eqnarray*}
|v\: m| &=& |v\: n|+|n\: m|=|v\: n|+|z\: n|-|z\: m|\\
&\leq& 2|v\: n|+|v\: z|-|m\: z|= |v\: x|+|v\: y|-|x\: y|+|v\: z|-|m\: z|\\
&=& \big(|v\: x|+|v\: y|+|v\: z|\big)-\big(|m\: x|+|m\: y|+|m\: z|\big)
\end{eqnarray*}
The first line uses Lemma~\ref{colinear} twice. The second line uses a triangle inequality, followed by formula~(\ref{leg}). The last line is a mere rewriting.

ii) We first prove the required inequality in the case when $y=y'$ and $z=z'$. By Lemma~\ref{colinear} we have
\[|m\: m'|=|m\: x'|-|m'\: x'|\leq |x\: x'|+|m\: x|-|m'\: x'|\]
and, similarly, 
\[|m\: m'|=|m'\: x|-|m\: x|\leq |x\: x'|+|m'\: x'|-|m\: x|.\]
Adding these two inequalities, we obtain the required $|m\: m'|\leq |x\: x'|$.

Next, we prove the required inequality in the case when $z=z'$. Let $n$ denote the median of $x,y',z$. By the previous step, we have $|n\: m'|\leq |x\: x'|$ and $|n\: m|\leq |y\: y'|$. It follows that $|m\: m'|\leq |x\: x'|+|y\: y'|$.

Finally, we prove the general case. Let $n$ denote the median of $x,y',z'$. Then $|n\: m'|\leq |x\: x'|$ by the first step, and $|n\: m|\leq |y\: y'|+|z\: z'|$ by the second step. We conclude that $|m\: m'|\leq |x\: x'|+|y\: y'|+|z\: z'|$, as required. \end{proof}

\begin{prop}\label{complete} The completion of a median space is a median space.
\end{prop}

\begin{proof} Let $X$ be a dense subspace of $\overline{X}$. Assuming that $X$ is a median space and $\overline{X}$ is complete, we need to show that $\overline{X}$ is a median space. 

Let $\overline{x}, \overline{y}, \overline{z}\in \overline{X}$. We show the existence of $\overline{m}\in\overline{X}$ such that $\overline{m}\in[\overline{x},\overline{y}]$, $\overline{m}\in[\overline{y},\overline{z}]$, $\overline{m}\in[\overline{z},\overline{x}]$. Indeed, let $x_k\to\overline{x}$, $y_k\to\overline{y}$, $z_k\to\overline{z}$ with $x_k,y_k,z_k\in X$, and let $m_k\in X$ denote the median of  $x_k,y_k,z_k$. By Lemma~\ref{ineq} ii), the sequence $(m_k)$ is a Cauchy sequence, so it converges to some $\overline{m}\in\overline{X}$. From $|x_k\: m_k|+|m_k\: y_k|=|x_k\: y_k|$ we get $|\overline{x}\;\overline{m}|+|\overline{m}\;\overline{y}|=|\overline{x}\;\overline{y}|$, i.e. $\overline{m}\in[\overline{x},\overline{y}]$, by letting $k\to\infty$. Similarly, $\overline{m}\in[\overline{y},\overline{z}]$ and $\overline{m}\in[\overline{z},\overline{x}]$.

Next, we argue that $\overline{m}\in\overline{X}$ is unique. Let $\overline{n}\in\overline{X}$ be such that $\overline{n}\in[\overline{x},\overline{y}]$, $\overline{n}\in[\overline{y},\overline{z}]$, $\overline{n}\in[\overline{z},\overline{x}]$, and let $n_k\to \overline{n}$ with $n_k\in X$. From Lemma~\ref{ineq} i) we have:
\[|n_k\: m_k|\leq \big(|n_k\: x_k|+|n_k\: y_k|+|n_k\: z_k|\big)-\big(|m_k\: x_k|+|m_k\: y_k|+|m_k\: z_k|\big)\] 
The right hand side converges to $0$ as $k\to\infty$; indeed, each parenthesized expression equals $\frac{1}{2}\big(|\overline{x}\;\overline{y}|+|\overline{y}\;\overline{z}|+|\overline{z}\;\overline{x}|\big)$ in the limit. Therefore $\overline{n}=\overline{m}$. \end{proof}


\subsection{Median spaces are negative definite} 
A metric space $X$ is said to be \emph{negative definite} if 
\[\sum\alpha_i\alpha_j |x_i\: x_j|\leq 0\] 
whenever $x_1, \ldots, x_n\in X$, and $\alpha_1, \ldots, \alpha_n\in\R$ with $\sum \alpha_i=0$. 

The geometry of negative definite metric spaces is related to the Hilbert space geometry by the following GNS construction. Assume that the metric space $X$ is negative definite. Let $V(X)$ be the real vector space on $X$, and let $V_0(X)$ consist of the vectors in $V(X)$ with zero coefficient sum. Then on $V_0(X)$ we can define an inner product as follows:
\[\Big\langle\sum\alpha_ix_i, \sum \beta_j y_j \Big\rangle=-\frac{1}{2}\sum\alpha_i\beta_j |x_i\:y_j|\] 
Note that we have $|x\: y|=\|x-y\|^2$.

Conversely, assume $X$ is a metric space such that there is an inner product space $V$ and a map $\gamma: X\To V$ with $|x\: y|=\|\gamma(x)-\gamma(y)\|^2$ for all $x,y\in X$. Then $X$ is negative definite. Indeed, for $\sum \alpha_i=0$ we can write
\[\sum\alpha_i\alpha_j  |x_i\: x_j|=-2\Big\|\sum\alpha_i\gamma(x_i)\Big\|^2\leq 0.\]

Before we handle the main proposition of this subsection, let us introduce the following definition. We say that $x,y,z,t\in X$ form a \emph{rectangle} if $y,t\in [x,z]$ and $x,z\in [y,t]$. In a rectangle, opposite sides have equal length: $|x\: y|=|z\:t|$, and $|y\:z|=|t\: x|$. 


\begin{prop}\label{negdef} Median spaces are negative definite.
\end{prop}
\begin{proof} Let $x_1, \ldots, x_n\in X$, and $\alpha_1, \ldots, \alpha_n\in\R$ with $\sum \alpha_i=0$. The median closure of a finite set is finite (Lemma~\ref{finite_median}), so we may assume that $\{x_1, \dots, x_n\}$ is median stable. Viewing $\{x_1,\ldots,x_n\}$ as a median space, we may furthermore assume that $X=\{x_1,\ldots,x_n\}$.

We argue by induction. Let $H$ be a maximal proper halfspace of $X$. 

For $x\in H^c$ let $p_x$ be a point in $H$ closest to $x$. We claim that $p_x\in [x,p]$ for all $p\in H$; in particular, $p_x$ is unique. Indeed, for $p\in H$ the median $m(p,p_x,x)$ is in $H$ and closer to $x$, unless $m(p,p_x,x)=p_x$, i.e., $p_x\in [x,p]$.

For every $x\in H^c$ we have that $H$ and $H^c$ are the only halfspaces that separate $x$ from $p_x$. Indeed, let $H'$ be a halfspace such that $p_x\in H'$ and $x\notin H'$. Pick $p\in H$; if $p\notin H'$ then $p_x\in [x,p]\subseteq H'^c$, contradiction. Thus $H\subseteq H'$. As $H'$ is proper ($x\notin H'$), we obtain $H=H'$ from the maximality of $H$. 

Next, we claim that $x,y,p_y, p_x$ forms a rectangle whenever $x,y\in H^c$. We know that $p_x\in [x,p_y]$ and $p_y\in [y,p_x]$. Assume $y\notin [x,p_y]$. Then there is a halfspace separating $[x,p_y]$ from $y$. As $p_x\in [x,p_y]$, the halfspace does not separate $x$ from $p_x$. On the other hand, the halfspace separates $y$ from $p_y$. This is a contradiction. Thus $y\in [x,p_y]$, and $x\in [y,p_x]$ by symmetry.

The consequence is that $|x\:p_x|$ is constant, say equal to $\delta$, for $x\in H^c$. 

Let $p:X\To H$ be the retraction that sends each $x\in X$ to the unique point in $p_x\in H$ that is closest to $x$. This extends the meaning of $p_x$ for $x\in H^c$ by setting $p_x=x$ for $x\in H$. We want to relate $|x\;y|$ and $|p_x\;p_y|$ for $x,y\in X$. There are several cases to consider:

$x\in H, y\in H$:\quad $|x\;y|=|p_x\;p_y|$ (trivially);

$x\notin H, y\notin H$:\quad $|x\;y|=|p_x\;p_y|$  (for $x, y, p_y, p_x$ forms a rectangle for $x,y\notin H$)

$x\notin H, y\in H$:\quad $|x\;y|=|p_x\;p_y|+\delta$ (for $p_x\in [x,y]$ gives $|x\;y|=|x\;p_x|+|p_x\;y|$)

$x\in H, y\notin H$:\quad $|x\;y|=|p_x\;p_y|+\delta$ (by symmetry with the previous case)

\noindent We now expand and replace as follows:
\begin{eqnarray*}
\sum \alpha_i\alpha_j|x_i\:x_j| &=&\sum_{x_i,x_j\in H} \alpha_i\alpha_j |x_i\:x_j| +\sum_{x_i\in H,x_j\notin H} \alpha_i\alpha_j |x_i\:x_j|+\\
&&\sum_{x_i\notin H,x_j\in H} \alpha_i\alpha_j |x_i\:x_j|+
\sum_{x_i,x_j\notin H} \alpha_i\alpha_j |x_i\:x_j|\\
&=& \sum_{x_i,x_j\in H} \alpha_i\alpha_j |p_{x_i}\:p_{x_j}| +\sum_{x_i\in H,x_j\notin H} \alpha_i\alpha_j \big(|p_{x_i}\:p_{x_j}|+\delta\big)+\\
&&\sum_{x_i\notin H,x_j\in H} \alpha_i\alpha_j \big(|p_{x_i}\:p_{x_j}|+\delta\big)+\sum_{x_i,x_j\notin H} \alpha_i\alpha_j |p_{x_i}\:p_{x_j}|\\
&=& \sum
\alpha_i\alpha_jd(p_{x_i},p_{x_j})+2 \delta\big(\sum_{x_i\notin H} \alpha_i\big)\big(\sum_{x_i\in H}\alpha_i\big)
\end{eqnarray*}
The first term, $\sum \alpha_i\alpha_jd(p_{x_i},p_{x_j})$, is non-positive by the induction hypothesis when applied to $H$. As $\sum \alpha_i=0$, the second term equals $-2\delta(\sum_{x_i\in H}\alpha_i)^2$. Thus $\sum \alpha_i\alpha_jd(x_i,x_j)\leq 0$. \end{proof}

\begin{rem}
A metric space $X$ is \emph{hypermetric} if $\sum t_it_j |x_i\: x_j|\leq 0$ whenever $x_1, \ldots, x_n\in X$, and integers $t_1, \ldots, t_n$ with $\sum t_i=1$. One easily shows that hypermetric spaces are negative definite. An obvious adaptation of the previous proof shows that median spaces are in fact hypermetric.
\end{rem}

\subsection{Group actions} Let $G$ act by isometries on a median space $X$ and consider the inner product space $V_0(X)$. The obvious linear action $g(\sum\alpha_ix_i):=\sum\alpha_i (gx_i)$ of $G$ on $V_0(X)$ is isometric, and for a fixed $v\in X$ we ``affinize'':
\begin{displaymath}
g*(\sum\alpha_ix_i)=g(\sum\alpha_ix_i)+(gv-v)
\end{displaymath}
The $*$ action extends to the Hilbert space completion of $V_0(X)$. More importantly, we have
$\|g*0\|=\|gv-v\|=\sqrt{|v\:gv|}$. We conclude:

\begin{prop}\label{FH_FM} 
If $G$ is a Kazhdan group, then every action of $G$ on a median space is bounded.
\end{prop}

\begin{prop}\label{PM_PH} 
If $G$ is a group acting properly on a median space, then $G$ is a Haagerup group.
\end{prop}

These two propositions generalize some of the statements that appear in Sections 2 and 3 (Theorems~\ref{kazhdan_FA}, \ref{kazhdan_FAR}, \ref{PQ-PH}).

\begin{conj}\label{conject} The converses of Proposition~\ref{FH_FM} and
Proposition~\ref{PM_PH} are true.
\end{conj}

\subsection{Verheul's work} As far as we know, the first thorough study of median geometry is that of Verheul \cite{Ver93}. Verheul's results seem to be relatively unknown; we only found out about them towards the end of this work. 

The fact that the completion of a median space is a median space (Proposition~\ref{complete}) appears as Corollary II.3.2 in \cite{Ver93}. Our proof seems to be different, and much more direct. This is partly due to the fact that Verheul works in a more general context, that of ``multimedian geometry''. Following the terminology of \cite{Ver93}, a metric space $X$ is \emph{modular} if, for each triple $x,y,z\in X$, the geodesic intervals $[x,y]$, $[y,z]$, $[z,x]$ have at least one common point. It is immediate to see that $X$ is a modular space if and only if \emph{Helly's property}
\begin{quote}
$C_1\cap\dots\cap C_n\neq \emptyset$ whenever $C_1,\dots,C_n$ are convex subsets that are mutually intersecting
\end{quote}
holds in $X$. Helly's property has been recently used in \cite{Fa08}.

One could also view Proposition~\ref{negdef} as being implicit in Verheul's work. Namely, the following statement is the core of Theorem V.2.4 in \cite{Ver93}: 

\begin{thm}[Verheul] Every median space isometrically embeds in some $L_1$ space.
\end{thm}

To obtain Proposition~\ref{negdef}, note that an $L_1$ space is negative definite; indeed, one need only integrate the fact that $\R$ is negative definite. Our proof of Proposition~\ref{negdef} is independent of Verheul's considerations, and it seems to us - once again - much simpler.


\section{Median graphs} 
Connected simplicial graphs come equipped with a path metric. One may thus consider median graphs. Immediate examples of median graphs are simplicial trees, the 1-skeleton of the square tiling of the plane, and the 1-skeleton of an $n$-dimensional cube. Note that median graphs are bipartite.

This short section describes the connection between spaces with walls, median graphs, and CAT(0) cube complexes.

\subsection{Spaces with walls} This notion was introduced by Haglund and Paulin \cite{HP98}. In our definition we insist on the presence of the trivial wall, which is needed for functorial reasons.
\begin{defn}\label{wallspacedef}
Let $X$ be a set. A \emph{wall} in $X$ is a partition of $X$ into 2 subsets called \emph{halfspaces}. We say that $X$ is a \emph{space with walls} if $X$ is endowed with a collection of walls, containing the trivial wall $\{\emptyset, X\}$, and so that any two distinct points are separated by a finite, non-zero
number of walls. 

A \emph{morphism of spaces with walls} is a map $f:X\To X'$ between spaces with walls with the property that $f^{-1}(A')$ is a halfspace of $X$ for each halfspace $A'$ of $X'$.
\end{defn}
Let $X$ be a space with walls. For $x\in X$, we let $\sigma_x$ denote the collection of halfspaces containing $x$. Then
\[|x\:y|_w=\frac{1}{2}|\sigma_x\triangle\sigma_y|\]
defines a metric on $X$, called the  \emph{wall metric}. Note that $|x\:y|_w$ equals the number of
walls separating $x$ and $y$.

A median graph $X$ is a space with walls, in the sense that the vertex set of $X$ together with the halfspaces afforded by the median structure satisfy Definition~\ref{wallspacedef}. One can show that the proper halfspaces in $X$ are all of the form $H_{xy}=\{z\in X^0:\: |z\:x|<|z\:y|\}$, where $xy$ is an edge. Consequently, the path metric and the wall metric coincide in a median graph.

\subsection{From spaces with walls to median graphs}\label{cubulation}
Every median graph can be viewed as a space with walls. Conversely, every space with walls admits a canonical embedding in a median graph in such a way that the wall structure is preserved. This canonical transformation of a space with walls into a median graph, an abstract version of a construction due to Sageev \cite{Sag95}, is explained in \cite{Nic04} and \cite{CN05}. We refer the reader to \cite{Nic04} for details on this brief presentation.

Out of a space with walls $X$ one constructs a connected median graph $\mathcal{C}^1(X)$. There is an embedding $X\hookrightarrow \mathcal{C}^1(X)$ which identifies $X$ as part of the vertex set of $\mathcal{C}^1(X)$. This embedding has the following properties:
\begin{itemize}
\item $X\hookrightarrow \mathcal{C}^1(X)$ is a morphism of spaces with walls; furthermore, the walls in $\mathcal{C}^1(X)$ identify bijectively with the walls of $X$ and, consequently, $X\hookrightarrow \mathcal{C}^1(X)$ is an isometric embedding when $X$ is equipped with the wall metric;

\item $X\hookrightarrow \mathcal{C}^1(X)$ is a ``dense'' embedding, in the sense that $\mathcal{C}^1(X)$ is the median closure of the image of $X$; in particular, if $X$ is already a median graph then
$X\hookrightarrow \mathcal{C}^1(X)$ is a median isomorphism, hence a graph isomorphism as well.
\end{itemize}
Also, $\mathcal{C}^1(\cdot)$ is functorial: a morphism of spaces with walls $f:X'\To X''$ extends uniquely to a median morphism $f_*:\mathcal{C}^1(X')\To \mathcal{C}^1(X'')$ (which need not be a graph morphism).

\subsection{From median graphs to CAT(0) cube complexes} 
The notation $\mathcal{C}^1(X)$ is meant to suggest that there is a natural complex $\mathcal{C}(X)$ having $\mathcal{C}^1(X)$ as its 1-skeleton. This complex $\mathcal{C}(X)$ is built from $\mathcal{C}^1(X)$ by ``filling in'' isometric copies of euclidean cubes, that is, by inductively adding an $(n+1)$-dimensional cube whenever its $n$-skeleton is present. One calls $\mathcal{C}(X)$ the \emph{cubulation} of $X$.

The crucial fact is that $\mathcal{C}(X)$ is a CAT(0) cube complex:

\begin{thm}[\cite{Rol98}, \cite{Che00}, \cite{Ger98}]\label{cubing_skeleton}
The $1$-skeleton of a CAT(0) cube complex is a median graph. Conversely, each
median graph is the $1$-skeleton of a CAT(0) cube complex, obtained by  ``filling in'' with higher euclidean cubes.
\end{thm}

This close relationship with CAT(0) cube complexes is one of the motivating facts about median graphs.

\subsection{Group actions}
Suppose $G$ acts on a space with walls $X$ by permuting the walls. Consequently, $G$ acts by isometries when $X$ is equipped with the wall metric. Let $\mathcal{H}$ denote the collection of halfspaces of $X$. Then $G$ acts linearly isometrically on $\ell_2(\mathcal{H})$ via $g\phi(H)=\phi(g^{-1}H)$.

Fix a basepoint $v\in X$ and consider the affine action:
\begin{displaymath}
g*\phi=g\phi+(\chi_{\sigma_{gv}}-\chi_{\sigma_{v}})
\end{displaymath}
Then
$\|g*0\|^2=\|\chi_{\sigma_{gv}}-\chi_{\sigma_{v}}\|^2=|\sigma_{gv}\triangle\sigma_v|=2|v\; gv|_w$.
We conclude:

\begin{prop}\label{FH-FW} 
Every action of a Kazhdan group on a space with walls is bounded.
\end{prop}

\begin{prop}\label{PH-PW} 
A group that admits a proper action on a space with walls is a Haagerup group.
\end{prop}

These propositions, implicit in \cite{NRe97}, are due to Haglund, Paulin and Valette.

The usefulness of Proposition~\ref{FH-FW} and Proposition~\ref{PH-PW} resides in the fact that spaces with walls can be read off in many geometric contexts. For instance, finitely generated Coxeter groups act properly on spaces with walls; there are at least three ways of interpreting their wall structure (see \cite{NRe03}).

The facts listed in Section~\ref{cubulation} imply that a group action on a space with walls $X$ extends uniquely to a group action on the median graph $\mathcal{C}^1(X)$. The embedding $X\hookrightarrow \mathcal{C}^1(X)$ is isometric, and a point $x\in X$ has the same orbit whether viewed in $X$ or $\mathcal{C}^1(X)$. It follows that the action on $X$ is bounded, respectively proper, if and only if the extended action on $\mathcal{C}^1(X)$ is bounded, respectively proper. The upshot is that Proposition~\ref{FH-FW}, respectively ~\ref{PH-PW}, are instances of Proposition~\ref{FH_FM}, respectively ~\ref{PM_PH}.

In light of Propositions~\ref{FH-FW} and ~\ref{PH-PW}, one is led to the following:

\begin{prob}\label{KW} i) Give examples of non-Kazhdan groups with the property that every action on a space with walls is bounded. ii) Give examples of Haagerup groups that do not admit proper actions on spaces with walls.
\end{prob}

Part ii) was recently solved by Haglund \cite{Hag07}. Part i) has an interesting reinterpretation, which we now explain.

\subsection{Sageev's theorem} 
Let $G$ be a finitely generated group. A subgroup $H\leq G$ is a \emph{codimension-1 subgroup} if the coset graph of $H$, i.e., the quotient of the Cayley graph of $G$ by $H$, has more than one end. Recall, the number of ends of a connected, locally finite graph $X$ is the supremum over the number of infinite components one obtains by excising arbitrary finite sets from $X$.

As observed by P. Scott, if $G$ splits over $H$ then $H$ is a codimension-1 subgroup of $G$.

Sageev \cite{Sag95} originally introduced the cubulation procedure in order to relate existence of codimension-1 subgroups with non-trivial actions on CAT(0) cube complexes. We now state a version of his theorem (see \cite{Ger98}, \cite{NRo98}, \cite{Rol98}). 

\begin{thm}\label{relative_ends}
Let $G$ be a finitely generated group. Then $G$ has no codimension-1 subgroups if and only if every action of $G$ on a median graph has bounded orbits.
\end{thm}

The above theorem, combined with Proposition~\ref{FH-FW} or Proposition~\ref{FH_FM}, shows that Kazhdan groups have no codimension-1 subgroups; this was first observed in \cite{NRo98}. Asking whether the converse holds amounts to Problem~\ref{KW} i).


\section{Appendix}
\subsection{Uniform convexity}\label{UC} Recall, a normed space $X$ is \emph{uniformly convex} if for
each $\e >0$ there is $\delta>0$ such that the following holds for all $x,y\in X$: 
\[
\|x\|\leq 1, \|y\|\leq 1, \|x-y\|\geq \e\;\Longrightarrow\;\Big\|\frac{1}{2}(x+y)\Big\|\leq 1-\delta
\]
This notion has a straightforward generalization to the metric context:

\begin{defn} A uniquely geodesic metric space $X$ is \emph{uniformly convex} if for
each $\e >0$ there is $\delta>0$ such that the following holds for all $x,y,z\in X$: 
\[
|x\;y|\geq \e \max\{|z\;x|, |z\;y|\}\;\Longrightarrow\;|z\;m_{xy}|\leq (1-\delta)\max\{|z\;x|, |z\;y|\}
\]
Here $m_{xy}$ denotes the midpoint of $x$ and $y$.
\end{defn}

$\cat$ metric spaces are uniformly convex. Indeed, a $\cat$ space $X$ is uniquely geodesic, and the CN inequality
\[|z\;m_{xy}|\leq\sqrt{\frac{|z\:x|^2+|z\:y|^2}{2}-\frac{|x\;y|^2}{4}}\]
holds in $X$. Uniform convexity with $\delta(\e)=\sqrt{1-\frac{1}{4}\e^2}$ follows immediately.

It is well-known that, in a uniformly convex Banach space or in a $\cat$ complete metric space, non-empty bounded sets have unique circumcenters. Recall, a point $x_*\in X$ is a \emph{circumcenter} for the non-empty bounded subset $A$ if $x_*$ minimizes the circumradius function $r(x)=\sup_{a\in A} |x\;a|$.

The next result unifies these two instances:

\begin{prop} In a uniformly convex, complete metric space, every non-empty bounded subset has a unique circumcenter.
\end{prop}

\begin{proof} Let $A$ be a non-empty bounded subset of the uniformly convex, complete metric space $X$. Set $r=\inf_{x\in X}\; r(x)$ with $r(x)=\sup_{a\in A} |x\;a|$ as above. The case when $A$ is a singleton is trivial, so we may assume $r>0$.

We show existence of circumcenters. Let $(x_n)$ be a sequence in $X$ with $r(x_n)\To r$. We claim that   $(x_n)$ is Cauchy. By completeness, it follows that $(x_n)$ converges to some $x_*\in X$. As the circumradius function is continuous, in fact $|r(x)-r(y)|\leq |x\;y|$, we obtain $r(x_*)=r$ and we conclude that $x_*$ is a circumcenter of $A$. 

The following claim implies the fact that $(x_n)$ is Cauchy:
\begin{quote} for each $\e>0$ there is $N_\e$ such that $|x_m\:x_n|<\e\max\{r(x_m),r(x_n)\}$ for all $m,n\geq N_\e$
\end{quote}
So let us argue the claim. Let $\e>0$. Pick $N_\e$ such that $r(x_n)\leq r/\sqrt{1-\delta}$ for all $n\geq N_\e$, where $\delta$ comes from the uniform convexity. Assume, by contradiction, that there are $m,n\geq N_\e$ so that $|x_m\:x_n|\geq\e\max\{r(x_m),r(x_n)\}$. Let $y_{m,n}$ denote the midpoint of $x_m$ and $x_n$. For each $a\in A$ we have $|a\;x_m|\leq r(x_m)$ and $|a\;x_n|\leq r(x_n)$, so the uniform convexity gives 
\[|a\; y_{m,n}|\leq (1-\delta)\max\{r(x_m),r(x_n)\}\leq r\sqrt{1-\delta}\]
which leads to the contradiction $r\leq r(y_{m,n})\leq r\sqrt{1-\delta}$.

We prove uniqueness of circumcenters. Assume $x_*$ and $x'_*$ are distinct circumcenters, and let $y$ be their midpoint. Pick $\e>0$ with $|x_*\:x'_*|\geq \e r$. For each $a\in A$ we have $|a\;x_*|\leq r(x_*)=r$, $|a\:x'_*|\leq r(x'_*)=r$, so $|a\;y|\leq (1-\delta)r$ for some $\delta=\delta(\e)>0$. Then $r\leq r(y)\leq (1-\delta)r$, a contradiction.
\end{proof}

\begin{rem}
1) In a proper metric space, the existence of circumcenters for non-empty bounded sets follows from the  observation that the circumradius function does not change by passing from a bounded set to its closure. 

2) Strict convexity, i.e., $|z\;m_{xy}|<\max\{|z\;x|, |z\;y|\}$ whenever $x\neq y$, suffices for uniqueness of circumcenters.
\end{rem}

Existence and uniqueness of circumcenters for non-empty bounded subsets has a standard consequence for isometric group actions:

\begin{cor} If an isometric group action on a uniformly convex, complete metric space has bounded orbits, then in fact the action has a (globally) fixed point.
\end{cor}

\subsection{Mazur-Ulam}\label{MU} We include, for the reader's convenience, a proof of the Mazur-Ulam theorem. This proof is inspired by \cite{Va03}. 

Recall that a map $f:X\To Y$ between real normed spaces is \emph{affine} if $f$ is linear up to a
translation, equivalently, $f$ satisfies $f((1-t)x+ty)=(1-t)f(x)+tf(y)$ for all $x,y\in X$ and $t\in [0,1]$. In order to show that a continuous function $f:X\To Y$ between real normed spaces is affine, it suffices to verify that $f(\frac{x+y}{2})=\frac{f(x)+f(y)}{2}$ for all $x,y\in X$.

\begin{thm}
Every bijective isometry between real normed spaces is affine.
\end{thm}

\begin{proof} 
Fix $x,y$ in a real normed space $X$. For every bijective isometry
$f$ defined on $X$, denote the possible ``affine defect'' by $\delta
(f)=\big\|f(\frac{x+y}{2})-\frac{f(x)+f(y)}{2}\big\|$. Note that we have a uniform bound on $\delta$:
\begin{displaymath}
\delta(f)\leq\frac{\|f(\frac{x+y}{2})-f(x)\|+\|f(\frac{x+y}{2})-f(y)\|}{2}=\frac{\|x-y\|}{2}
\end{displaymath}
For every bijective isometry $f$ defined on $X$ consider $f^\dagger:=f^{-1}\rho f$,
where $\rho$ is the reflection in $\frac{f(x)+f(y)}{2}$ in the target
space of $f$. Then:
\begin{eqnarray*}
\delta(f^\dagger)&=&\Big\|f^{-1}\Big(f(x)+f(y)-f\Big(\frac{x+y}{2}\Big)\Big)-\frac{y+x}{2}\Big\|\\
&=&\Big\|f(x)+f(y)-2f\Big(\frac{x+y}{2}\Big)\Big\|=2\delta (f)
\end{eqnarray*}
It follows that every bijective isometry defined on $X$ has vanishing affine defect.
\end{proof}


\end{document}